\documentclass[11pt]{amsart}
\usepackage{amsthm, amsfonts, amssymb, amsmath, bbold, tkz-berge, caption}\usetikzlibrary{decorations.pathreplacing, graphs}
\usepackage{raleway}

\newcommand{\G}{\Gamma}
\newcommand{\C}{\mathbb{C}}

\newcommand{\F}{\mathbb{F}}

\newcommand{\la}{\langle}
\newcommand{\ra}{\rangle}

\newcommand{\ct}{\mathbb{1}}

\newtheorem{thm}{Theorem}[section]

\theoremstyle{definition}

\begin{document}
\title[Independence number of regular graphs of matrix rings]{On the independence number of regular graphs of matrix rings}

\date{\today}
\author{Bogdan Nica}

\begin{abstract}  
Consider a graph on the non-singular matrices over a finite field, in which two distinct non-singular matrices are joined by an edge whenever their sum is singular. We prove an upper bound for the independence number of this graph. As a consequence, we obtain a lower bound for its chromatic number that significantly improves a previous result of Tomon.
\end{abstract}

\address{\newline Department of Mathematical Sciences \newline Indiana University Indianapolis}
\email{bnica@iu.edu}

\maketitle

\section{Introduction}

Let $M_n(\F_q)$ denote the ring of $n\times n$ matrices over $\F_q$, the finite field with $q$ elements. The \emph{regular graph} over $M_n(\F_q)$ has vertex set $GL_n(\F_q)$, the multiplicative group of non-singular matrices in $M_n(\F_q)$. Two distinct matrices $a,b\in GL_n(\F_q)$ are joined by an edge whenever $a+b$ is singular, that is, $\det(a+b)=0$. In what follows, we denote the regular graph over $M_n(\F_q)$ by $\G_n(q)$. The case $n=1$ being rather trivial, we assume $n\geq 2$ in what follows.

Regular graphs, and some of their subgraphs, have been investigated in several papers. A result due to Akbari, Jamaali, and Seyed Fakhari \cite{AJ} says that, when $q$ is odd, the clique number of the regular graph $\G_n(q)$ satisfies
\[\omega(\G_n(q))\leq \sum_{k=0}^n k! \binom{n}{k}^2.\]
Interestingly, this bound is independent of $q$.

The chromatic number of the regular graph $\G_n(q)$ was studied by Tomon \cite{T}. For $q$ odd, he showed that
\[
\chi(\G_n(q))\geq (q/4)^{\lfloor n/2 \rfloor},
\]
though he believed that the above lower bound could be improved. In this note we deduce a much better lower bound for $\chi(\G_n(q))$, on the order of $q^{n-1}$. Additionally, we remove the parity restriction on $q$; in fact, for even $q$ one could get an even better lower bound. Tomon also showed an upper bound, $\chi(\G_n(q))=O(q^{n^2-n})$, for $q$ odd. It would be interesting to further narrow this wide gap between the lower and upper bounds.

For bounds on the chromatic number for certain subgraphs of regular graphs, see the works of Akbari, Aryapoor, and Jamaali \cite{AA}, respectively Bardestani and Mallahi-Karai \cite{B}.

Our main result is the following upper bound on the independence number of the regular graph $\G_n(q)$.

\begin{thm}\label{main}
We have
\[\alpha(\G_n(q))\leq q^{n^2-n+1}.\]
\end{thm}

A lower bound on the chromatic number is then an immediate consequence, as $\alpha(\G_n(q))\cdot \chi(\G_n(q))\geq |GL_n(\F_q)|$. We have
\[|GL_n(\F_q)|=(q^n-1)(q^n-q)\dots(q^n-q^{n-1})=c_n(q)\cdot q^{n^2}\]
where
\begin{align}\label{eq: coeff}
c_n(q)=(1-q^{-1})\dots (1-q^{-n}).
\end{align}
Clearly, $0<c_n(q)<1$ for each $n\geq 1$; in fact
\[1-q^{-1}-q^{-2}<c_n(q)\leq 1-q^{-1}.\]
The lower bound for $c_n(q)$, the most useful one for our purposes, follows from the stronger inequality $1-q^{-1}-q^{-2}+q^{-(n+1)}\leq c_n(q)$; this, in turn, can be easily checked by induction.

We therefore deduce the following consequence. 

\begin{thm}\label{thm: chrom}
We have
\[\chi(\G_n(q))\geq (1-q^{-1}-q^{-2})q^{n-1}.\]
\end{thm}

The regular graphs under consideration are instances of a general graph-theoretic construction over rings. Let $R$ be a ring, and let $Z(R)$ denote the set of zero-divisors in $R$. The complement, $R\setminus Z(R)$, is the set of regular elements of $R$. The \emph{regular graph} of $R$ has vertex set $R\setminus Z(R)$, and two distinct elements $a,b\in R\setminus Z(R)$ are joined by an edge whenever $a+b\in Z(R)$. A related graph is the \emph{total graph} of $R$, having vertex set $R$ and the same adjacency law. Thus, the total graph of $R$ contains the regular graph of $R$ as an induced subgraph. Total and regular graphs over rings were introduced by Anderson and Badawi \cite{AB}, in the commutative context, and then by Dol\v{z}an and Oblak \cite{DO} in general. For an overview of the study of graphs over rings, we refer to the recent monograph \cite{A+}.

Returning to the ring of interest, $M_n(\F_q)$, let us denote the corresponding total graph by $T_n(q)$. To clarify, the vertex set of $T_n(q)$ is $M_n(\F_q)$, and two distinct matrices $a,b\in M_n(\F_q)$ are joined by an edge whenever $\det(a+b)=0$. As we explain below, the total graph $T_n(q)$ plays a key role in our approach.

We prove Theorem~\ref{main} by spectral methods. The most direct attack would be to invoke Hoffman's ratio bound. Recall, this says the following: if $X$ is a regular graph of degree $d$ on $v$ vertices, then 
\begin{align}\label{eq: H}
\alpha(X)\leq\frac{v}{1-d/\theta_{\min}}
\end{align}
where $\theta_{\min}$ is the smallest adjacency eigenvalue of $X$; see Haemers \cite{H} for a nice account. The regular graph $\G_n(q)$ may be viewed as the Cayley graph of the general linear group $GL_n(\F_q)$ with respect to the symmetric subset $S'=\{s\in GL_n(\F_q): \det (I_n+s)=0\}$; if $q$ is even, then the identity $I_n$ has to be removed from $S'$. This shows, firstly, that the regular graph $\G_n(q)$ is regular in the graph-theoretic sense. It is, however, rather non-trivial to evaluate the degree of $\G_n(q)$--that is to say, the size of $S'$. Secondly, the Cayley graph viewpoint offers a way to compute the adjacency spectrum of the regular graph $\G_n(q)$, and the conjugation-invariance of the set $S'$ is convenient to that end. Alas, this approach requires the character theory of $GL_n(\F_q)$.

We circumvent these difficulties by following a similar strategy in the total graph $T_n(q)$. The total graph, a supergraph of  
the regular graph, turns out to be more amenable to spectral computations. The reason is that its underlying set, $M_n(\F_q)$, is an abelian group under addition, so its character theory is much simpler. We will prove that, in fact, the following holds.
\begin{thm}\label{thm: tot}
We have
\[\alpha(T_n(q))\leq q^{n^2-n+1}.\]
\end{thm}
Theorem~\ref{main} clearly follows from Theorem~\ref{thm: tot}.

\section{Proof of Theorem~\ref{thm: tot}}
\subsection{Sum-graphs} Let $G$ be a finite abelian group, whose operation is written additively, and let $S$ be a subset of $G$. The \emph{sum-graph} of $G$ with respect to $S$ has vertex set $G$, and two distinct elements $g,h\in G$ are connected whenever $g+h\in S$. 

The sum-graph construction is a variation on the Cayley graph construction. Two differences should be highlighted. Firstly, the sum-graph construction yields a a simple graph--undirected and loopless--no matter what $S$ is; for the Cayley graph construction, one needs $S$ to be symmetric. Secondly, a Cayley graph is regular, but a sum-graph need not be. In fact, the sum-graph of $G$ with respect to $S$ is nearly $|S|$-regular: the vertex degrees are $|S|$ or $|S|-1$. Vertices having degree $|S|-1$ correspond to those $g\in G$ which satisfying $2g\in S$. Two opposite cases are worth singling out:

\begin{itemize}
\item[$\bullet$] the `binary' case, when $G$ is an elementary $2$-group, meaning that $2g=0$ for all $g\in G$; in this case, the sum-graph of $G$ with respect to $S$ is the Cayley graph of $G$ with respect to $S\setminus\{0\}$;
\item[$\bullet$] the `odd' case, when $|G|$ is odd; then there are $|S|$ vertices of degree $|S|-1$, and the remaining $|G|-|S|$ vertices have degree $|S|$.
\end{itemize}

The total graph $T_n(q)$ is the sum-graph of the abelian group $M_n(\F_q)$ with respect to the subset of singular matrices 
\[S=\{s\in M_n(\F_q): \det(s)=0\}.\]

When $q$ is even, we are in the `binary' case; the total graph $T_n(q)$ is regular of degree $|S|-1$ in this case. When $q$ is odd, we are in the `odd' case; the total graph $T_n(q)$ is no longer regular, in fact the singular matrices have degree $|S|-1$, while the non-singular matrices have degree $|S|$.

\subsection{Laplacian spectrum of sum-graphs} The sum-graph structure of the total graph $T_n(q)$ unlocks the computation of its Laplacian spectrum. Indeed, the following general result shows how to determine the Laplacian eigenvalues of a sum-graph by using the characters of the underlying abelian group.

\begin{thm}\label{spec}
Let $G$ be a finite abelian group, and let $S\subseteq G$. Then the Laplacian eigenvalues of the sum-graph of $G$ with respect to $S$ are given as follows: $|S|-\sum_{s\in S} \phi(s)$ for each real character $\phi$, respectively $|S|\pm |\sum_{s\in S} \phi(s)|$ for each conjugate pair of non-real characters $\{\phi,\overline{\phi}\}$.
\end{thm}
 
\begin{proof} By definition, the Laplacian operator $L$ acts on a function $\phi:G\to \C$ as follows: for each vertex $g\in G$, we have
\begin{align*}
(L\phi)(g)=\deg(g) \phi(g)-\sum_{h:\: h\sim g} \phi(h)
\end{align*}
where $\deg(g)$ denotes the degree of the vertex $g$, and the latter sum is taken over the vertices adjacent to $g$. The degree of a vertex $g$ is either $|S|$, and then the neighbors of $g$ are $\{s-g:s\in S\}$, or $|S|-1$, in which case the neighbors of $g$ are $\{s-g:s\in S\}\setminus\{g\}$. At any rate, we have
\begin{align*}
(L\phi)(g)=|S|\phi(g)-\sum_{s\in S} \phi(s-g).
\end{align*}
Now if $\phi:G\to\C$ is a character, then
\begin{align*}
\sum_{s\in S} \phi(s-g)=\bigg(\sum_{s\in S} \phi(s)\bigg)\overline{\phi}(g).
\end{align*}
By using the shorthand $\la \phi, S\ra=\sum_{s\in S} \phi(s)$ we may write, concisely,
\begin{align*}
L\phi=|S|\phi-\la \phi, S\ra \overline{\phi}.
\end{align*}

It follows that, with respect to a suitable indexing of the basis provided by the characters, $L$ takes a block-diagonal form. Each real character $\phi$ contributes a diagonal entry $|S|-\la \phi, S\ra$, and each conjugate pair of non-real characters $\{\phi,\overline{\phi}\}$ contributes a $2\times 2$ matrix
\begin{align*}
\renewcommand\arraystretch{1.4}
\begin{pmatrix}
|S| & -\la \phi, S\ra\\
-\overline{\la \phi, S\ra} & |S|
\end{pmatrix}.
\end{align*}
The eigenvalues of the above matrix are $|S|\pm |\la \phi, S\ra|$. The spectrum of $L$ is now easily read off. 
\end{proof}

Theorem~\ref{spec} is an analogue of the well-known recipe for computing spectra of Cayley graphs over abelian groups. In the regular context, passing between adjacency and Laplacian spectra is trivial. While there is no apparent recipe for computing the adjacency spectrum of a sum-graph, let us point out that, thanks to near-regularity, the adjacency eigenvalues can be well-approximated by using the Laplacian eigenvalues.

If we remove the word `distinct' from the definition of a sum-graph, then the resulting graph is no longer simple--there may be loops sprouting at certain vertices. In the literature, one usually encounters the `loopy' version of sum-graphs. It is not hard to see that computing the Laplacian spectrum of a (loopless) sum-graph is equivalent to computing the adjacency eigenvalues of the `loopy' sum-graph. The adjacency spectrum for `loopy' sum-graphs has been known for some time:  see \cite[Prop.1]{Li}, \cite[Thm.2.1]{D+}; compare also the weaker version in \cite[Sec.4]{Chu}. Theorem~\ref{spec} is therefore equivalent to a known result. Yet, it appears to be formally novel--to the best of our knowledge, it has not been explicitly stated, proved, or used in this form before.

There is at least one real character of $G$, the trivial character $\ct$; fittingly, it yields the trivial Laplacian eigenvalue $\lambda=0$ for any sum-graph over $G$. In the `binary' case, each non-trivial character is real; in the `odd' case, no non-trivial character is real. 

\subsection{Laplacian eigenvalues of total graphs} We turn to applying Theorem~\ref{spec} to the total graph $T_n(q)$. The set $S$ is the complement of $GL_n(\F_q)$ in $M_n(\F_q)$, so 
\[|S|=q^{n^2}-|GL_n(\F_q)|=q^{n^2}-c_n(q)\cdot q^{n^2}\]
with $c_n(q)$ as in \eqref{eq: coeff}.

The characters of the additive group $M_n(\F_q)$ can be described by means of a fixed non-trivial additive character $\psi$ of $\F_q$. For each matrix $u\in M_n(\F_q)$, define $\phi_u: M_n(\F_q)\to \C$ by 
\[\phi_u(x)=\psi\bigg(\sum_{i,j=1}^n u_{ij}\;x_{ij}\bigg).\]
Then $\phi_u$ is a character of $M_n(\F_q)$, and each character is of this form for a unique $u\in M_n(\F_q)$. Formally, $u\mapsto \phi_u$ is an isomorphism between the additive group $M_n(\F_q)$ and its dual.

We also need the following evaluation, due to Li and Hu \cite[Thm.2.1]{LH}: if $u\in M_n(\F_q)$ is a non-zero matrix of rank $r$, then 
\begin{align}
\sum_{a\in GL_n(\F_q)} \phi_u(a)=(-1)^r \: q^{n(n-1)/2}\: \prod_{i=1}^{n-r} (q^i-1).
\end{align}
As $\sum_{a\in M_n(\F_q)} \phi_u(a)=0$, we deduce that
 \begin{align}
\sum_{s\in S} \phi_u(s)=-(-1)^r \: q^{n(n-1)/2}\: \prod_{i=1}^{n-r} (q^i-1).
\end{align}

Putting everything together, we deduce the following.

\begin{thm}\label{thm: LS}
The non-trivial distinct Laplacian eigenvalues of the total graph $T_n(q)$ are as follows:
\begin{itemize}
\item[$\bullet$] for $q$ even, 
\[|S|+(-1)^r \: q^{n(n-1)/2}\: \prod_{i=1}^{n-r} (q^i-1), \qquad r=1,\dots,n;\]
\item[$\bullet$] for $q$ odd, 
\[|S|\pm \: q^{n(n-1)/2}\: \prod_{i=1}^{n-r} (q^i-1), \qquad r=1,\dots,n.\]
\end{itemize}
\end{thm}

Theorem~\ref{thm: LS} implies that the number of distinct Laplacian eigenvalues of $T_n(q)$ is $n+1$ when $q$ is even, respectively $2n+1$ when $q$ is odd. It is a general fact that the number of distinct Laplacian eigenvalues of a graph is greater than its diameter. Here, however, we note that the diameter of $T_n(q)$ is very small--namely, $2$--for all $n\geq 2$. For if $a,b\in M_n(\F_q)$ are two distinct vertices in $T_n(q)$, then consider a vertex $c\in M_n(\F_q)$ such that $a+c$ has zero first row, while $b+c$ has zero second row;  if $c$ is distinct from $a$ and $b$, then $c$ is adjacent to both $a$ and $b$, whereas if $c=a$ or $c=b$, then $a$ and $b$ are already adjacent.

\subsection{A Hoffman-type bound} Next, we wish to invoke the following generalization of Hoffman's ratio bound \eqref{eq: H}: if $X$ is a graph of minimal degree $\delta$ on $v$ vertices, then 
\begin{align}\label{eq: HT}
\alpha(X)\leq v \Big(1-\frac{\delta}{\lambda_{\max}}\Big)
\end{align}
where $\lambda_{\max}$ is the largest Laplacian eigenvalue of $X$. The Hoffman-type bound \eqref{eq: HT} was first proved by van Dam and  Haemers \cite[Lem.3.1]{vDH}, and then rediscovered several times over, e.g., \cite[Cor.3.6]{GN}. In the case of regular graphs, \eqref{eq: HT} turns into the Hoffman bound \eqref{eq: H}. See also \cite{N} for a refinement of \eqref{eq: HT}.

The final step towards the proof of Theorem~\ref{thm: tot} is applying the Hoffman-type bound \eqref{eq: HT} to the total graph $T_n(q)$. We have $v=q^{n^2}$ and $\delta\geq |S|-1$, while $\lambda_{\max}$ can be read off from Theorem~\ref{thm: LS}. For the sake of uniformity, we simply estimate
\begin{align*}
\lambda_{\max}&\leq |S|+q^{n(n-1)/2}\: \prod_{i=1}^{n-1} (q^i-1)\\
&< |S|+q^{n(n-1)/2}\prod_{i=1}^{n-1} q^i=|S|+q^{n^2-n}.
\end{align*}
(We could do better in the case when $q$ is even.) Therefore
\[1-\frac{\delta}{\lambda_{\max}}<1-\frac{|S|-1}{|S|+q^{n^2-n}}=\frac{q^{n^2-n}+1}{|S|+q^{n^2-n}}.\]
But $|S|=(1-c_n(q))\cdot q^{n^2}>q^{-1}\cdot q^{n^2}=q^{n^2-1}$, and so
\[\frac{q^{n^2-n}+1}{|S|+q^{n^2-n}}<\frac{q^{n^2-n}+1}{q^{n^2-1}+q^{n^2-n}}\leq \frac{1}{q^{n-1}}.\]
We conclude, as claimed, that
\[\alpha(T_n(q))\leq q^{n^2}\cdot\frac{1}{q^{n-1}}=q^{n^2-n+1}.\]

\end{document}